\documentclass{ifacconf}

\usepackage{graphicx}      
\usepackage{epstopdf}
\DeclareGraphicsExtensions{.pdf,.jpeg,.png,.eps}

\usepackage{natbib}        
\begin{document}
\begin{frontmatter}

\title{Efficient particle continuation model predictive control} 

\thanks[footnoteinfo]{Accepted to the 16th IFAC Workshop on Control Applications
of Optimization (CAO'2015), Garmisch-Partenkirchen, Germany, October 6--9, 2015.}

\author[First]{Andrew Knyazev} 
\author[Second]{Alexander Malyshev}

\address[First]{Mitsubishi Electric Research Labs (MERL)
201 Broadway, 8th floor, Cambridge, MA 02139, USA. (e-mail: knyazev@merl.com), \hfill
({http://www.merl.com/people/knyazev}).}
\address[Second]{Mitsubishi Electric Research Labs (MERL)
201 Broadway, 8th floor, Cambridge, MA 02139, USA. (e-mail: malyshev@merl.com)}
\begin{abstract}                
Continuation model predictive control (MPC), introduced by T. Ohtsuka in 2004,
uses Krylov-Newton approaches to solve MPC optimization and is suitable for nonlinear
and minimum time problems. We suggest particle continuation MPC in the case,
where the system dynamics or constraints can discretely change on-line.
We propose an algorithm for on-line controller implementation of continuation MPC
for ensembles of predictions corresponding to various anticipated changes
and demonstrate its numerical effectiveness for a test minimum time problem arriving
to a destination. Simultaneous on-line particle computation of ensembles of controls,
for several dynamically changing system dynamics, allows choosing the optimal destination
on-line and adapt it as needed.
\end{abstract}


\end{frontmatter}

\section{Introduction}

Model predictive control (MPC) is a popular control approach, which efficiently
treats constraints on state and control variables. Solid introduction into MPC
is found in \cite{CaBo:04} and \cite{GrPa:11}, industrial applications are discussed in
\cite{QiBa:03}, numerical aspects of MPC are surveyed in \cite{DiFeHa:09} and
\cite{WaBo:10}.

\cite{Oht:04} has developed an on-line numerical method for nonlinear MPC,
based on the so-called Newton-Krylov method; see \cite{KnKe:04} for other
applications of the Newton-Krylov method. \cite{KnFuMa:15} provide
an efficient preconditioner for the Newton-Krylov method and extend the MPC model,
treated in \cite{Oht:04}, to the minimum-time problem.

In the present note, we demonstrate how the preconditioned Newton-Krylov method
by \cite{KnFuMa:15} works in cases, where the system dynamics or constraints can
discretely change on-line. The problem with discrete switches requires simultaneous
solution to several finite-horizon predictions, which can be done independently
on parallel processors. The computer time can be also reduced by using the same
preconditioner for all finite-horizon control problems, if the discrete changes
of the system dynamics or constraints do not lead to large norms of the residual
mapping $F$, introduced in Section 2.
If the norm $\|F\|$ is not small enough after a discrete switch,
the method may be even ruined during execution. Such a behavior is inherited
from the Newton method, which converges only for sufficiently good initial guess.
The finite-horizon predictions with $\|F\|$ larger than a suitable tolerance
may be discarded or refined.

The necessity to compute several close trajectories also appears in the so called
particle control problems; see e.g. \cite{OhOh:11}.

The rest of the note is organized as follows. Section 2 presents a framework of
the receding horizon prediction problem and obtains its solution in the form of
a nonlinear equation. Section 3 discusses how this nonlinear equation is efficiently solved
for an ensemble of sufficiently close trajectories issued from the current state. 
Section 4 describes a test example and presents formulas for computer implementation.
Section 5 shows numerical results.

\section{Receding horizon prediction}

In this section, we introduce the receding horizon problem for an MPC model
considered in~\cite{KnFuMa:15}. The finite horizon is an interval $[t,t+T]$,
where $T$ may depend on $t$.
The control input $u(\tau)$ and parameter vector~$p$ are determined so as to
minimize the performance index
\[
J = \phi(x(t+T),p)+\int_t^{t+T}L(\tau,x(\tau),u(\tau),p)d\tau
\]
subject to the state dynamics
\begin{equation}\label{e1}
\frac{dx}{d\tau}=f(\tau,x(\tau),u(\tau),p),
\end{equation}
and the equality constraints for the state $x$ and the control~$u$
\begin{equation}\label{e2}
C(\tau,x(\tau),u(\tau),p) = 0,
\end{equation}
\begin{equation}\label{e3}
\psi(x(t+T),p) = 0.
\end{equation}
The initial value condition $x(\tau)|_{\tau=t}$ for equation (\ref{e1})
is the state vector $x(t)$ of the dynamic system. The control vector $u=u(\tau)|_{\tau=t}$,
solving the problem over the receding horizon, is used afterwards as an input
to control the dynamic system at time~$t$.

Let us discretize the continuous control problem stated above on a time grid
$\tau_i$, $i=0,1,\ldots,N$, obtained by partitioning the horizon $[t,t+T]$ into
$N$ subintervals of size $\Delta\tau_i=\tau_{i+1}-\tau_i$. The vector
functions $x(\tau)$ and $u(\tau)$ are replaced by their values $x_i$ and $u_i$
at the grid points $\tau_i$. The integral in the performance index $J$ is approximated
by the rectangular quadrature rule. Equation (\ref{e1}) is integrated by the explicit
Euler method. The discretized optimal control problem is as follows:
\[
\min_{u_i,p}\left[
\phi(x_N,p) + \sum_{i=0}^{N-1}L(\tau_i,x_i,u_i,p)\Delta\tau_i\right],
\]
subject to
\begin{equation}\label{e4}
\quad x_{i+1} = x_{i}+f(\tau_i,x_i,u_i,p)\Delta\tau_i,\quad i = 0,1,\ldots,N-1,
\end{equation}
\begin{equation}\label{e5}
C(\tau_i,x_i,u_i,p) = 0,\quad  i = 0,1,\ldots,N-1,
\end{equation}
\begin{equation}\label{e6}
\psi(x_N,p) = 0.
\end{equation}

The necessary optimality conditions for the discretized finite horizon problem
are obtained by means of the discrete Lagrangian function
\begin{eqnarray*}
&&\mathcal{L}(X,U)=\phi(x_N,p)+\sum_{i=0}^{N-1}
L(\tau_i,x_i,u_i,p)\Delta\tau_i\\
&&+\,\lambda_0^T[x(t)-x_0]+\sum_{i=0}^{N-1}\lambda_{i+1}^T
[x_i-x_{i+1}+\Phi_i(\tau_i,x_i,u_i,p)\Delta\tau_i]\\
&&+\sum_{i=0}^{N-1}\mu_i^TC(\tau_i,x_i,u_i,p)\Delta\tau_i+\nu^T\psi(x_N,p),
\end{eqnarray*}
where $X = [x_i\; \lambda_i]^T$, $i=0,1,\ldots,N$, and
$U = [u_i\; \mu_i\; \nu\; p]^T$, $i=0,1,\ldots,N-1$.
Here, $\lambda$ is the costate vector and $\mu$ is the Lagrange multiplier vector 
associated with constraint~(\ref{e5}). The terminal constraint (\ref{e6})
is relaxed by the aid of the Lagrange multiplier $\nu$. 
The necessary optimality conditions are the KKT stationarity conditions:
$\mathcal{L}_{\lambda_i}=0$, $\mathcal{L}_{x_i}=0$, $i=0,1,\ldots,N$,
$\mathcal{L}_{u_j}=0$, $\mathcal{L}_{\mu_j}=0$, $i=0,1,\ldots,N-1$,
$\mathcal{L}_{\nu_k}=0$, $\mathcal{L}_{p_l}=0$.

For convenience in the subsequent formulations, we introduce the Hamiltonian function
\begin{eqnarray*}
\lefteqn{H(t,x,\lambda,u,\mu,p) = L(t,x,u,p)}\hspace*{8em}&& \\
&&{}+\lambda^T f(t,x,u,p)+\mu^T C(t,x,u,p).
\end{eqnarray*}

The number of unknowns in the KKT conditions can be reduced by eliminating
the states $x_i$ and costates $\lambda_i$:

(1) Starting from the current measured or estimated state $x_0$, compute
$x_i$, $i=0,1\ldots,N-1$, by the forward recursion
\[
x_{i+1} = x_i + f(\tau_i,x_i,u_i,p)\Delta\tau_i.
\]
Then starting from
\[
\lambda_N=\frac{\partial\phi^T}{\partial x}(x_N,p)+
 \frac{\partial\psi^T}{\partial x}(x_N,p)\nu
\]
compute the costates $\lambda_i$, $i=N\!-\!1,\ldots,0$, by the backward recursion
\[
\lambda_i=\lambda_{i+1}+\frac{\partial H^T}{\partial x}
(\tau_i,x_i,\lambda_{i+1},u_i,\mu_i,p)\Delta\tau_i.
\]
(2) Combine the control input $u$, the Lagrange multiplier
$\mu$, the Lagrange multiplier $\nu$, and the parameter $p$, all in one vector
\[
U(t)=[u_0^T,\ldots,u_{N-1}^T,\mu_0^T,\ldots,\mu_{N-1}^T,\nu^T,p^T]^T. 
\]
Calculate the mapping $F[U,x,t]$, using just obtained values $x_i$ and $\lambda_i$, as

\begin{eqnarray*}
\lefteqn{F[U,x,t]}\\
&&\hspace*{-2em}=\left[\begin{array}{c}\begin{array}{c}
\frac{\partial H^T}{\partial u}(\tau_0,x_0,\lambda_{1},u_0,\mu_0,p)\Delta\tau_0\\
\vdots\\\frac{\partial H^T}{\partial u}(\tau_i,x_i,\lambda_{i+1},u_i,\mu_i,p)\Delta\tau_i\\
\vdots\\\frac{\partial H^T}{\partial u}(\tau_{N-1},x_{N-1},\lambda_{N},u_{N-1},
\mu_{N-1},p)\Delta\tau_{N-1}\end{array}\\\;\\
\begin{array}{c}C(\tau_0,x_0,u_0,p)\Delta\tau_0\\
\vdots\\C(\tau_i,x_i,u_i,p)\Delta\tau_i\\\vdots\\
C(\tau_{N-1},x_{N-1},u_{N-1},p)\Delta\tau_{N-1}\end{array}\\\;\\
\psi(x_N,p)\\[2ex]
\begin{array}{c}\frac{\partial\phi^T}{\partial p}(x_N,p)+
\frac{\partial\psi^T}{\partial p}(x_N,p)\nu\\
+\sum_{i=0}^{N-1}\frac{\partial H^T}{\partial p}(\tau_i,x_i,
\lambda_{i+1},u_i,\mu_i,p)\Delta\tau_i\end{array}
\end{array}\right].
\end{eqnarray*}
The vector argument $x$ in $F[U,x,t]$ denotes the initial vector $x_0$ in the forward recursion.

The equation with respect to the unknown vector $U(t)$
\[
 F[U(t),x(t),t]=0
\]
gives the required necessary optimality conditions that are solved
on the controller in real time.

\section{Numerical algorithm}

The controlled dynamic system is sampled on a discrete time grid $t_j$, $j=0,1,\ldots$.
The sampled values of the state and parameters are $x_j=x(t_j)$ and $U_j=U(t_j)$.
Choosing a small $h$, which is usually much less than $\Delta t_j=t_{j+1}-t_j$ and $\Delta\tau_k$,
we introduce the operator
\[
a_j(V)=(F[U_{j-1}+hV,x_j,t_j]-F[U_{j-1},x_j,t_j])/h.
\]
The discrete equation $F[U_j,x_j,t_j]=0$ is then equivalent to the operator equation
$a_j(\Delta U_j/h)=b_j/h$, where $\Delta U_j=U_j-U_{j-1}$, $b_j=-F[U_{j-1},x_j,t_j]$.
This operator equation allows us to compute $U_j$, if $U_{j-1}$ is known.

Let us denote the $k$-th column of the $m\times m$ identity matrix by $e_k$,
where $m$ is the dimension of the vector $U$, and form an $m\times m$ matrix $A_j$
with the columns $A_je_k=a_j(e_k)$, $k=1,\ldots,m$.
The matrix $A_j$ is an $O(h)$ approximation of the Jacobian matrix
$F_U[U_{j-1},x_j,t_j]$, which is symmetric.

We suppose that $U_0$ is an approximate solution to the equation $F[U_0,x_0,t_0]=0$ and
omit discussion of methods for computing $U_0$. The first block entry of $U_0$ is
taken as the control $u_0$ at the state $x_0$. The next state $x_1=x(t_1)$
is either sensor estimated or computed by the formula $x_1=x_0+f(t_0,x_0,u_0)\Delta t_0$;
cf. (\ref{e4}).

At the time $t_j$, $j>1$, we have the state $x_j$ and the vector $U_{j-1}$ evaluated
at the previous time $t_{j-1}$. We solve the following operator
equation with respect to $V$:
\begin{equation}\label{e7}
 a_j(V)=b_j/h.
\end{equation}
Then we set $\Delta U_j=hV$, $U_j=U_{j-1}+\Delta U_j$ and choose the first block
component of $U_j$ as the control $u_j$. The next system state $x_{j+1}=x(t_{j+1})$
is either sensor estimated or computed by the formula $x_{j+1}=x_j+f(t_j,x_j,u_j)\Delta t_j$.

An approximate solution to (\ref{e7}) can be found by computing the matrix $A_j$
and then solving the system of linear equations $A_j\Delta U_j=b_j$ by the Gaussian
elimination with pivoting.
A more efficient way is solving (\ref{e7}) by the GMRES method,
where the operator $a_j(V)$ is used instead of $A_j$
as in \cite{Oht:04} and \cite{KnFuMa:15}.

Convergence of GMRES can be accelerated by preconditioning.
A matrix $M$ that is close to the matrix $A$ and such that computing $M^{-1}r$
for an arbitrary vector $r$ is relatively easy, is referred to as a preconditioner.
The preconditioning for the system of linear equations $Ax=b$ with the preconditioner
$M$ formally replaces the original system $Ax=b$ with the equivalent preconditioned
linear system $M^{-1}Ax=M^{-1}b$. When the condition number $\|M^{-1}A\|\|A^{-1}M\|$
is sufficiently small, convergence of the preconditioned GMRES is fast.
The vector $z=M^{-1}r$ is often computed via back-substitutions as $z=U^{-1}(L^{-1}r)$,
where $L$ and $U$ are the triangular factors in the LU factorization $M=LU$ computed
by the Gaussian elimination.

\cite{KnFuMa:15} compute the matrix $A_j$ exactly for some time instances $t_j$
and use it as a preconditioner for GMRES at a number of subsequent time
instances $t_j$, $t_{j+1}$, \ldots, $t_{j+n_{prec}}$.

In the present paper, we suggest to apply the above described numerical method
in more general situations, where the dynamic system and/or constraints depend on
a discrete parameter, or switch, with few values $1$, \ldots, $q$.
In other words, there are $q$ functions $f^{(k)}(t,x,u,p)$ and $q$ mappings
$F^{(k)}(U,x,t)$, $k=1,2,\ldots,q$, and at each time instance $t_j$ we must select $k$
such that the performance index is minimized.

The numerical method is modified as follows. At the time $t_j$, $j>1$, the system state
is given by $x_j$, and we have the vector $U_{j-1}$ evaluated at the previous time $t_{j-1}$.
Since there are $q$ mappings $F^{(k)}$, we must solve $q$ operator equations with respect to $V$:
\begin{equation}\label{e8}
 a_j^{(k)}(V^{(k)})=b_j^{(k)}/h.
\end{equation}
Solutions $V^{(k)}$ with low precision can be discarded or refined.
In the set of admissible solutions $V^{(k)}$, we select the solution $V^{(k_j)}$ with
minimum performance index. 
Then we set $\Delta U_j=hV^{(k_j)}$, $U_j=U_{j-1}+\Delta U_j$ and choose the first block
component of $U_j$ as the control $u_j$. The next system state $x_{j+1}=x(t_{j+1})$
is either sensor estimated or computed by the formula
$x_{j+1}=x_j+f^{(k_j)}(t_j,x_j,u_j)\Delta t_j$.

Equations (\ref{e8}) are solved by GMRES independently. However, a preconditioner
for all $k$ can be a single matrix $A_j$, which is evaluated for some suitable $k$.

The vector $U_{j-1}$ satisfies $F^{(k_{j-1})}(U_{j-1},x_{j-1},t_{j-1})\approx0$.
When each of the $q$ equations $F^{(k)}(U_j^{(k)},x_{j},t_{j})=0$ is solved by
the Newton method, it may succeed only if all the residuals
$F^{(k)}(U_{j-1},x_{j},t_{j})$ are sufficiently small, or at least some of them.
Thus, the changes by the discrete switches should not be
too radical, which is the main limitation of the modified method.

\section{Test problem}

We consider a minimum-time problem on the two-dimensional plane from a state
$(x_0,y_0)$ to a state $(x_f,y_f)$ with inequality constraints.
The system dynamics is governed by the system of differential equations
\begin{equation}\label{e10}
\frac{d}{dt}\left[\begin{array}{c}x\\y\end{array}\right]=
\left[\begin{array}{c}(Ax+B)\cos u\\(Ax+B)\sin u\end{array}\right].
\end{equation}
The control variable $u$ is subject to an inequality constraint:
$u$ stays within the band $c_{u}-r_{u}\leq u\leq c_{u}+r_{u}$.
Following \cite{Oht:04} we introduce a slack variable $u_s$ and replace
the inequality constraint by the equality
\[
C(u,u_d)=(u-c_{u})^2+u_s^2-r_{u}^2=0. 
\]
The state is forced to pass through the point $(x_f,y_f)$ at time $t=t_f$
by imposing two terminal constraints
\[
\psi(x,y,p)=\left[\begin{array}{c}x-x_f\\y-y_f\end{array}\right]=0. 
\]
The objective is to minimize the performance index
\[
J=\phi(p)+\int_{t_0}^{t_f}L(x,y,u,u_s,p)dt',
\]
where
\[
\phi(p)=p=t_f-t_0,\quad L(x,y,u,u_s,p)=-w_su_s.
\]
The term $\phi(p)$ is responsible for the shortest time to destination,
and the function $L$ serves to stabilize the slack variable $u_s$.

For convenience, we change the time variable $t$ within the horizon
by the new time $\tau=(t-t_0)/(t_f-t_0)$, which runs over the interval $[0,1]$.

The corresponding discretized finite-horizon problem on a uniform grid $\tau_i$
uses the following data structures and computations:
\begin{itemize}
\item $\tau_i=i\Delta\tau$, where $i=0,1,\ldots,N$, and $\Delta\tau=1/N$;
\item the participating variables are the state $\left[\begin{array}{c}
x_i\\y_i\end{array}\right]$, the costate $\left[\begin{array}{c}
\lambda_{1,i}\\\lambda_{2,i}\end{array}\right]$, the control
$\left[\begin{array}{c}u_{i}\\u_{si}\end{array}\right]$, the Lagrange multipliers
$\mu_i$ and $\left[\begin{array}{c}\nu_{1}\\\nu_{2}\end{array}\right]$;
\item the system dynamics is governed by the equations
\[
\left\{\begin{array}{l} x_{i+1}=x_i+\Delta\tau p(Ax_i+B)\cos u_{i},\\
\,y_{i+1}=y_i+\Delta\tau p(Ax_i+B)\sin u_{i},\end{array}\right.
\]
where $i=0,1,\ldots,N-1$;
\item the costate is computed by the backward recursion ($\lambda_{1,N}=\nu_1$,
$\lambda_{2,N}=\nu_2$)
\[
\left\{\begin{array}{l} \lambda_{1,i}=\lambda_{1,i+1}
-\Delta\tau pA(\cos u_i \lambda_{1,i+1}+\sin u_i\lambda_{2,i+1}),\\
\lambda_{2,i} = \lambda_{2,i+1}, \\
\end{array}\right.
\]
where $i=N-1,N-2,\ldots,0$;
\item the nonlinear equation $F(U,x_0,t_0)=0$, where
\begin{eqnarray*}
\lefteqn{U=[u_0,\ldots,u_{N-1},u_{s,0},\ldots,u_{s,N-1},}\hspace*{8em}\\
&&\mu_0,\ldots,\mu_{N-1},\nu_1,\nu_2,p],
\end{eqnarray*}
has the following rows from the top to bottom:
\[
\left\{\begin{array}{l}       
\Delta\tau[ p(Ax_i+B)\left(-\sin u_i\lambda_{1,i+1}+\cos u_i\lambda_{2,i+1}\right)\\
\hspace*{12em}{}+2\left(u_i-c_{u}\right)\mu_i] = 0\end{array}\hspace*{2em}\right.
\]
\[
\left\{\;\;\Delta\tau p\left[2\mu_iu_{si}-w_sp\right] = 0\hspace*{12em}\right.
\]
\[
\left\{\;\;\Delta\tau p\left[(u_i-c_{u})^{2}+u_{si}^2-r_{u}^2\right]=0
\right.\hspace*{8em}
\]
\[
\left\{\;\begin{array}{l}x_N-x_f=0\\y_N-y_f=0\end{array}\right.\hspace{17.5em}
\]
\[
\left\{\begin{array}{l}\Delta\tau \{\sum\limits^{N-1}_{i=0}(Ax_i+B)
(\cos u_i\lambda_{1,i+1}+\sin u_i\lambda_{2,i+1})\\
\hspace{12em}{}-w_su_{si}\}+1 = 0.\end{array}
\hspace*{4em}\right.
\]
\end{itemize}

\section{Numerical results}

We set $q=3$ and consider $q$ cases of the system dynamics simultaneously.
The cases are determined by the three pairs of constants:
$(A_1,B_1)=(0.97,1.)$, $(A_2,B_2)=(0.9,1.05)$, $(A_3,B_3)=(1.1,0.9)$.
Other parameters are the same in all cases: the end points of the computed trajectory are
$(x_0,y_0)=(0,0)$ and $(x_f,y_f)=(1,1)$; the constants in the inequality
constraint for the control are $c_u=0.8$ and $r_u=0.2$; $w_s=0.005$.

The number of grid points on the horizon is $N=20$,
the time step of the dynamic system is $\Delta t=1/200$,
the numerical differentiation step is $h=10^{-8}$.

The value of $U$ at time $t_0$ is approximated by the MATLAB
function \texttt{fsolve} with a special initial guess.

We use the GMRES method without restarts implemented in MATLAB.
The number of GMRES iterations does not exceed $30$, and the
absolute tolerance of the GMRES iterations equals $10^{-5}$.

We apply a simple preconditioning strategy as follows. The
exact Jacobian $F_U$ is computed periodically at time instances
with the period $0.2$. Then the LU factorization of
the Jacobian is used as the preconditioner until the next time
when it is recomputed.

\begin{figure}[ht]
\begin{center}
\includegraphics[width=\columnwidth]{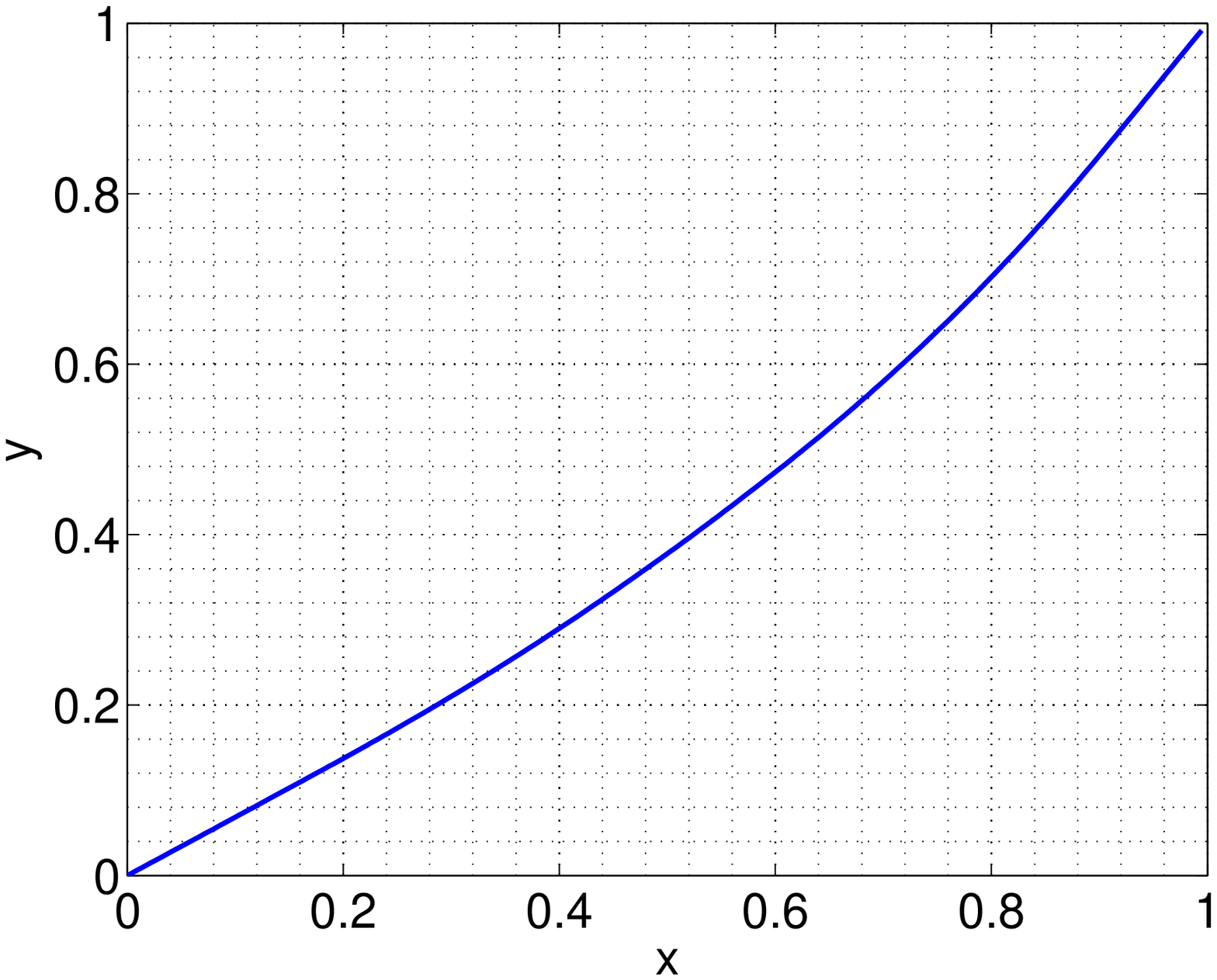}
\caption{The computed trajectory $(x,y)$.}
\label{fig1}
\end{center}
\end{figure}

\begin{figure}[ht]
\begin{center}
\includegraphics[width=\columnwidth]{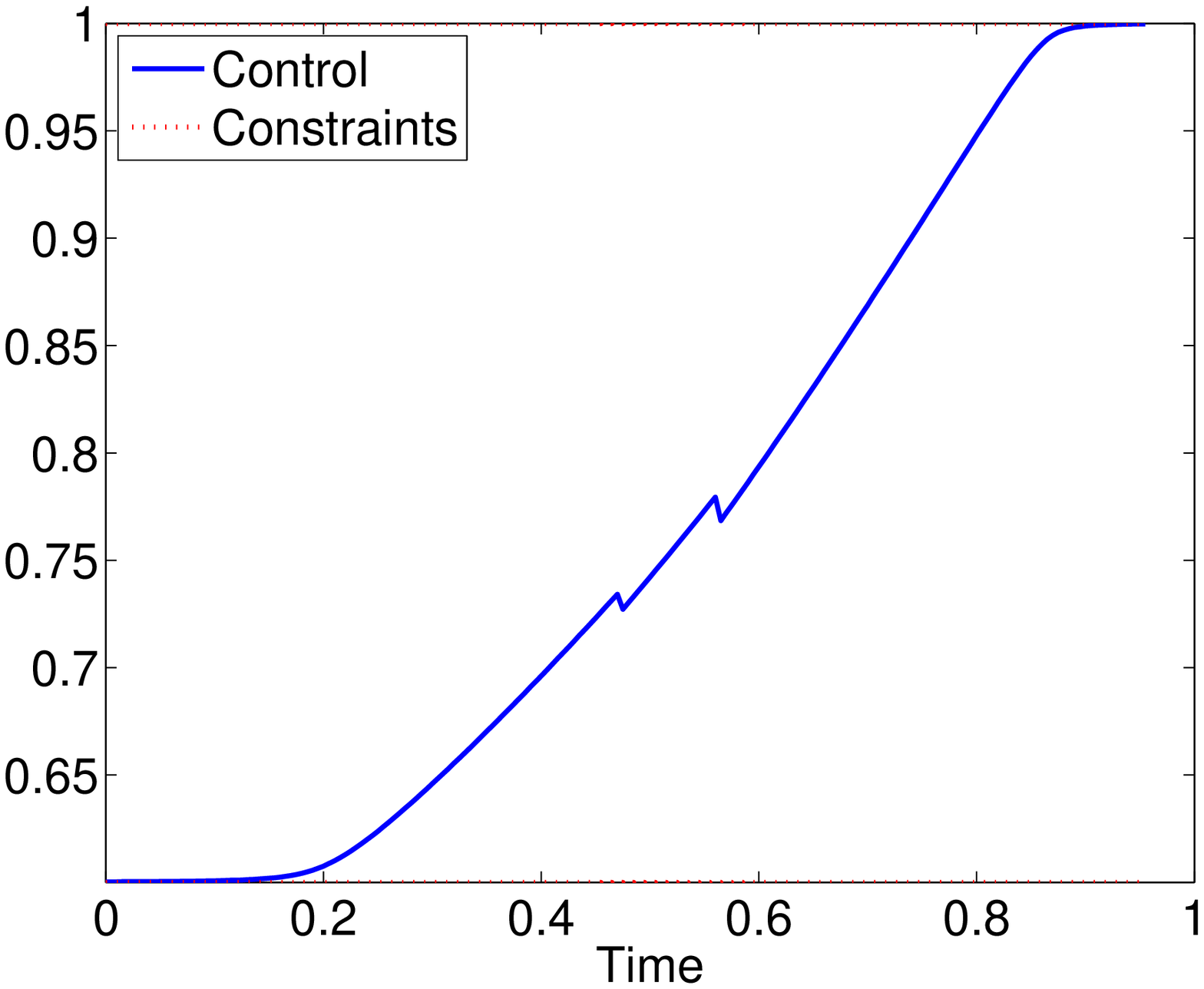}
\caption{The computed control $u(t)$.}
\label{fig2}
\end{center}
\end{figure}

\begin{figure}[ht]
\begin{center}
\includegraphics[width=\columnwidth]{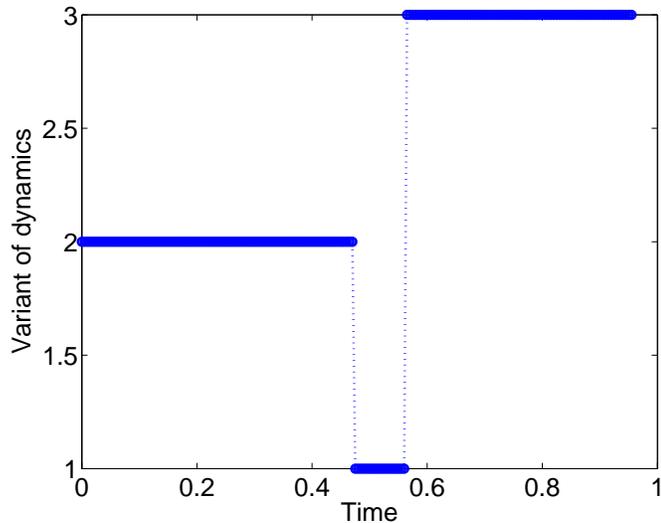}
\caption{The number $k$ of the chosen pair $(A_k,B_k)$.}
\label{fig3}
\end{center}
\end{figure}

\begin{figure}[ht]
\begin{center}
\includegraphics[width=\columnwidth]{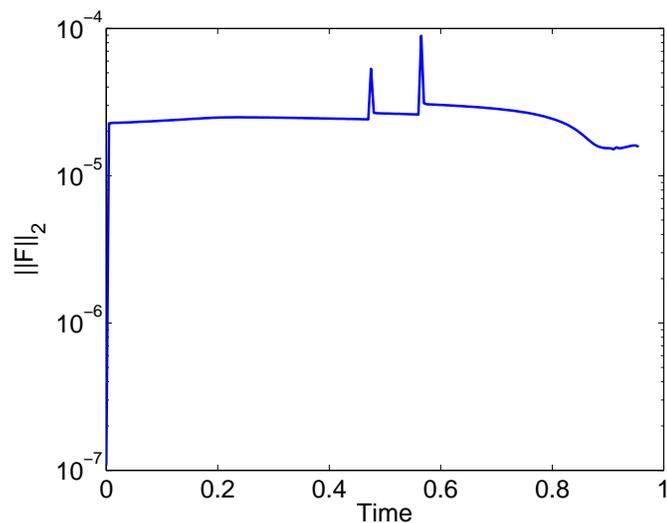}
\caption{The residual norm $\|F\|_2$.}
\label{fig4}
\end{center}
\end{figure}

Figure~\ref{fig1} displays the computed trajectory for the test problem with
automatic switches in the system dynamics. Time to destination along this trajectory
at the initial point $(x_0,y_0)$ is $0.974$. Figure~\ref{fig2} plots the control $u(t)$
computed by our method. Figure~\ref{fig3} shows how the system dynamics switch
between $(A_k,B_k)$ along the trajectory: first 95 steps are executed with $k=2$,
then 18 steps with $k=1$, and last 79 steps with $k=3$.
The residual norm $\|F\|_2$ is shown in Figure~\ref{fig4}.

We observe in Figure~\ref{fig2} that during a switch to another pair of $(A_k,B_k)$,
the control $u(t)$ undergoes abrupt change. The corresponding deterioration
of the norm $\|F\|_2$ at these points is seen in Figure~\ref{fig4}.
The method continues to work while the deterioration is sufficiently small.
A robust implementation of our method should verify $\|F\|_2$ when switching
between discrete parameters.

\begin{figure}[ht]
\begin{center}
\includegraphics[width=\columnwidth]{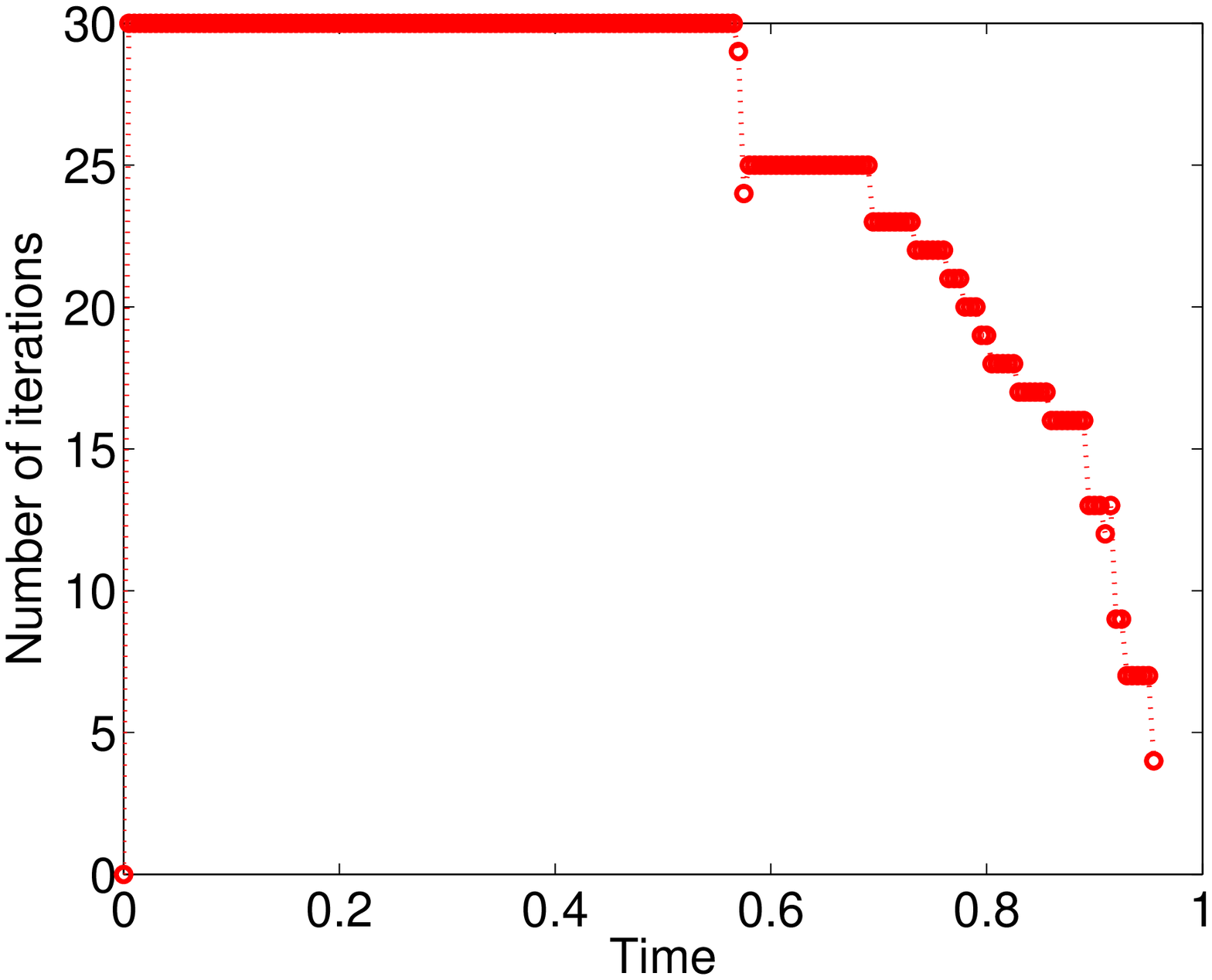}
\caption{The number of iterations for GMRES without preconditioning.}
\label{fig5}
\end{center}
\end{figure}

\begin{figure}[ht]
\begin{center}
\includegraphics[width=\columnwidth]{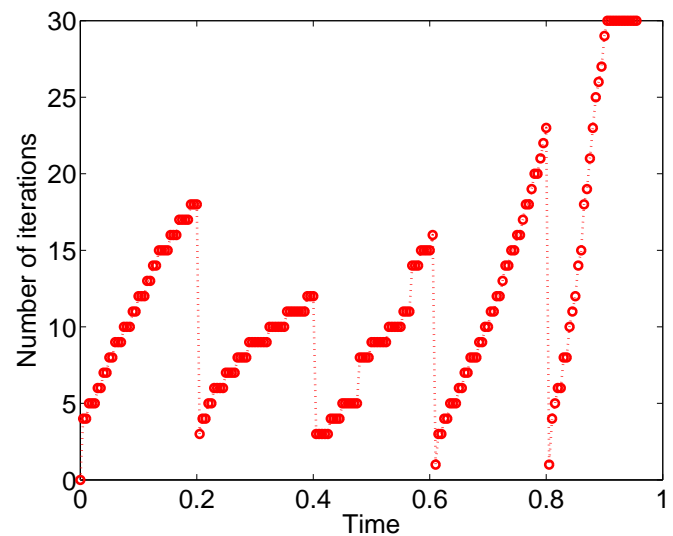}
\caption{The number of iterations for GMRES with preconditioning.}
\label{fig6}
\end{center}
\end{figure}

Figures~\ref{fig5} and \ref{fig6} show the number of GMRES iterations in the
non-preconditioned and preconditioned variants, respectively.
The variant without preconditioning uses 4920 iterations in total.
The variant with preconditioning uses 2216 iterations in total.

\section{Conclusion}

The numerical method, developed for nonlinear MPC problems in \cite{KnFuMa:15},
can be used in cases, when an ensemble of near solutions for the finite-horizon
prediction have to be computed simultaneously.
The reduction of computing time may be achieved by using parallel processors for
each prediction and/or by using a single preconditioner for the whole ensemble.

\end{document}